\documentclass[12pt]{article}
\usepackage{mathrsfs}
\usepackage{amsfonts}
\usepackage{amssymb,bbm}  
\usepackage{amsmath,amscd,graphicx,amsthm}
\usepackage{subfigure}
\usepackage{subeqnarray}
\usepackage{cases}
\usepackage{color}
\usepackage{paralist}  
\usepackage[colorlinks,urlcolor=black,anchorcolor=black,linkcolor=black,citecolor=blue]{hyperref}
\usepackage{amsmath}
\usepackage[numbers,sort&compress]{natbib}
\usepackage{indentfirst}
\makeatletter \@addtoreset{equation}{section}

       \newtheorem{De}{\bf Definition}[section]
       \newtheorem{lem}{\bf Lemma}[section]
       
       \newtheorem{expl}{\bf Example}[section]
       \newtheorem{assp}{\bf Assumption}

       \newtheorem{thm}{\bf Theorem}[section]

\newcommand{\R}{\mathbb{R}}

\def\be{\begin{equation}}     \def\ee{\end{equation}}
\def\bea{\begin{eqnarray}}    \def\eea{\end{eqnarray}}
\def\beaa{\begin{eqnarray*}}  \def\eeaa{\end{eqnarray*}}

\begin{document}
\title{The delay feedback control for the McKean-Vlasov stochastic differential equations with common noise}

\author{Xing Chen\thanks{School of Mathematics and Statistics,
Northeast Normal University, Changchun, Jilin, 130024, China. }
 \and Xiaoyue Li\thanks{School of Mathematical Sciences, Tiangong University, Tianjin, 300387, China. Research of this author  was  supported by the National Natural Science Foundation of China (No. 12371402, 11971096), the National Key R$\&$D Program of China (2020YFA0714102), the Natural Science Foundation of Jilin Province, China (No. YDZJ202101ZYTS154), and the Fundamental Research Funds for the Central Universities, China (No. 2412021ZD013).}
\and Chenggui Yuan\thanks{Department of Mathematics, Swansea University, Bay Campus, SA1 8EN, UK.}
}

\date{}

\maketitle
\begin{abstract}
 Since response lags are essential in the feedback loops and are required by most physical systems, it is more appropriate to stabilize McKean-Vlasov stochastic differential equations (MV-SDEs) with common noise through the implementation of delay feedback control mechanisms.  The aim of this paper is to design delay feedback control functions of the system state such that the controlled system to be boundedness in infinite horizon and further exponentially stable in the mean square. The designed controller, which depends only on the system state is easier to implement than that
   in \cite{WHGY2022} which was designed to depend  on both system state and measure. The existence and uniqueness of the global solution of the controlled system is proved. The It\^o formula with respect to both state and measure  is derived. The proposed delay feedback control strategies are rendered viable for effective stabilization of MV-SDEs with common noise. Furthermore, the moment Lyapunov exponent, which is intricately linked to the time delays, is meticulously estimated.

\end{abstract}

\section{Introduction}\label{sec:intr}
The McKean-Vlasov stochastic differential equations (MV-SDEs), the limiting equations of a particle system in mean-field particle systems as the number of particles tends to infinity, have received much attention (see, for example, \cite{Mc1966,S1991,W2018,DQ2021}). 
If the particles have a common source of randomness $W^0$, the limiting equations of a particle in mean-field particle systems becomes the MV-SDEs with common noise
\begin{align}\label{1}
\mathrm{d} X(t)
     &=f(X(t) ,\mathcal{L}^1( X(t)))
     \mathrm{d} t
+g( X(t) ,\mathcal{L}^1( X(t))) \mathrm{d} W_t
\nonumber\\&\quad
+g^0( X(t), \mathcal{L}^1( X(t))) \mathrm{d} W_t^0,
\end{align}
where $\mathcal{L}^1( X(t))$ represents the conditional distribution of $X(t)$ given the common noise source $W^0$, which have been applied to various fields including finance, physics, biology, mean-field games, machine learning, among others \cite{BSW2023,BLM2023,BC2023,CLZ2023,KT,KX1999,LL2023,LLW2022,W2018}.
The fundamental theory of the MV-SDEs with common noise \eqref{1} is well established. 
One can refer to  \cite{CD2018,KNRS2022} for the well-posedness and conditional propagation of chaos, \cite{HSS2021} for the weak well-posedness, \cite{M,BW2024} for the long-time behavior.

In numerous disciplines, including finance, physics, and biology, the automatic control constitutes a pivotal concern, with considerable attention given to the stabilization analysis \cite{CZS2009,DLM2012,M2015,YMMH2015,AK2014,LMMY20,MLH2008,LXZ2007,WHGY2022}.
For this purpose,
it is traditional to design a feedback control term, denoted as
$u(t,X(t))$,  such that the controlled system
 becomes stable. However, the requisite response lags for the majority of physical systems are indispensable in feedback loops \cite{RMT2015}.  Consequently, it is not only logical but also pertinent to design a control term $u(t,X(t-\tau))$ such that the controlled system
becomes stable, where $\tau$ signifies the temporal lag between the state observation and the system's response to the control. Although the stabilization by delay feedback control for stochastic differential equations has been discussed by several authors (see, \cite{YHLM2022,Z2019,LMMY20,MLH2008}, to name a few), there is so far little on the delay feedback control for the MV-SDEs with common noise.

The main aim of this paper is to design the delay feedback control  such that the controlled MV-SDE with common noise
\begin{align}\label{6.1*}
\mathrm{d} X(t)
&=\big(f(X(t) ,\mathcal{L}^1(X(t)))+u(t,X(t-\tau))\big)
     \mathrm{d} t
+g(X(t),\mathcal{L}^1(X(t))) \mathrm{d} W_t
\nonumber\\&\quad
+g^0(X(t), \mathcal{L}^1(X(t))) \mathrm{d} W_t^0,
\end{align}
becomes stable.
 Mathematically speaking, this paper
  employs the It\^o formula with respect to both state and measure and the Lyapunov functional method. It is worthy noting that  the introduction of common noise and the response lag increases the complexity of theoretical analysis. The analysis on the conditional distribution is rather technical. The main contributions of this paper are highlighted as follows.
\begin{itemize}
\item The existence and uniqueness of the global solution of the controlled system   is proved. The It\^o formula with respect to both state and measure of the control system is investigated;

\item Delay feedback control mechanisms are meticulously designed to ensure that the controlled MV-SDEs with common noise exhibit asymptotic boundedness and mean-square exponential stability, meanwhile, explicit bounds on the permissible delay $\tau$ are derived, and the impact of $\tau$ on the Lyapunov exponents is elucidated;
\item The control term   designed here only depends on the system state which is easier to implement than that designed in \cite{WHGY2022} for MV-SDEs.
\end{itemize}

The subsequent sections of this paper are organized as follows: Section 2 presents some necessary preliminaries. Section 3 shows the existence and uniqueness of solution of the controlled system, and derives the It\^o formula with respect to both state and measure of the controlled system. Section 4 is dedicated to the boundedness control for the MV-SDEs with common noise. Section 5 deals with the stabilization problem of the controlled MV-SDEs with common noise and  provides an example to illustrate the result on stabilization. Section 6 concludes this paper by summarizing the contributions.

\section{Preliminaries}

Throughout this paper we use the following notations. Denote by $|\cdot|$ the Euclidean norm in $\mathbb{R}^n$. If $A$ is a matrix or vector, let its transpose be $A^{\prime}$ and its trace norm be $|A|=\sqrt{\operatorname{tr}(A^{\prime}A)}$. If $A$ is constant, $\lceil A\rceil$ represents the rounding upwards of $A$. 
Let $(\Omega^0, \mathcal{F}^0, \mathbb{P}^0)$ and $(\Omega^1, \mathcal{F}^1,\mathbb{P}^1)$ be two complete probability spaces endowed with complete and right-continuous filtration $(\mathcal{F}_t^0)_{t \geq 0}$ and $(\mathcal{F}_t^1)_{t \geq 0}$, respectively. Here $(W^0_t)_{t \geq 0}$ and $(W_t)_{t \geq 0}$  are $d$-dimensional Brownian motions on  $\left(\Omega^0, \mathcal{F}^0, \mathbb{P}^0\right)$ and $\left(\Omega^1, \mathcal{F}^1, \mathbb{P}^1\right)$, respectively.  
 Define a product space $(\Omega, \mathcal{F}, \mathbb{P})$, where
$\Omega=\Omega^0 \times \Omega^1,$ {$(\mathcal{F},\mathbb{P})$} is the completion of $(\mathcal{F}^0\otimes \mathcal{F}^1,\mathbb{P}^0 \otimes \mathbb{P}^1)$, and $(\mathcal{F}_t)_{t\geq 0}$ is the complete and right-continuous {augmentation} of $(\mathcal{F}^0_t\otimes\mathcal{F}^1_t)_{t\geq 0}$. Denote the element of $\Omega$  by $\omega=(\omega^0,\omega^1)$   {where $\omega^0\in \Omega^0$ and   $\omega^1\in \Omega^1$}.  Let $\mathbb{E}^1$ be the expectation under $\mathbb{P}^1$. For $i=1,2$, let $({\Omega}^1_i,\mathcal{F}^1_i,\mathbb{P}^1_i)$ be the copies of $({\Omega}^1,\mathcal{F}^1,\mathbb{P}^1)$, $\mathbb{E}^1_i$ be the expectation under $\mathbb{P}^1_i$ and $\Omega_i=\Omega^0 \times {\Omega}^1_i$ be the copies of $\Omega$. Denoted by $\mathcal{P}(\mathbb{R}^d)$ the space of
all probability measures over $\mathbb{R}^d$ equipped with the weak topology. For any $k>0$, let
$$
\mathcal{P}_k(\mathbb{R}^d):=\left\{\mu \in \mathcal{P}(\mathbb{R}^d): {\int_{\mathbb{R}^d} |x|^k\mu(\mathrm{d}x) }<\infty\right\},
$$
and the  Wasserstein distance $\mathbb{W}_k$
$$
\mathbb{W}_k(\mu, \nu):=\inf_{\pi \in \Pi(\mu, \nu)}\left(\int_{\mathbb{R}^d \times \mathbb{R}^d}|x-y|^k \pi(\mathrm{d} x, \mathrm{~d} y)\right)^{\frac{1}{k\vee 1}},
$$
for any $\mu, \nu \in \mathcal{P}_k(\mathbb{R}^d)$, $\Pi(\mu, \nu)$ is the set of all couplings of $\mu$ and $\nu$. $\mathcal{P}_k(\mathbb{R}^d)$ is a Polish space under the  Wasserstein distance $\mathbb{W}_k$. 
Let $\delta_0$ be the Dirac measure at point $0\in\mathbb{R}^d$.

 Given a random variable $\xi$ on $\Omega$ equipped with the $\sigma$-algebra $\mathcal{F}^0\otimes \mathcal{F}^1$, for any $\omega^0\in\Omega^0$, $\xi(\omega^0,\cdot):\Omega^1\ni\omega^1\rightarrow X(\omega^0,\omega^1)$ is a random variable on $\Omega^1$ and $\mathcal{L}(\xi(\omega^0,\cdot))$ is the law of $\xi(\omega^0,\cdot)$. For  a random variable $\xi$ on $(\Omega, \mathcal{F}, \mathbb{P})$, $\mathcal{L}^1(\xi):=\mathcal{L}(\xi(\omega^0,\cdot))$ forms a random variable from $\Omega^0$ to $\mathcal{{P}}(\mathbb{R}^d)$, which provides a conditional law of $\xi$ given $\mathcal{F}^0$(see, \cite[Vol-II, Lemma 2.4]{CD2018}).  For $k\geq1$, define $\mathbb{L}_k(\mathbb{R}^d)$ by the family of all $\mathbb{R}^d$-valued random variables $\xi$ 
 such that $\mathbb{E}|\xi|^k<\infty$.

It is worth noting that we do not distinguish a random variable $X$ on $(\Omega^0, \mathcal{F}^0, \mathbb{P}^0)$ from its natural extension $\tilde{X}:(\omega^0,\omega^1)\mapsto X(\omega^0)$ on $(\Omega, \mathcal{F}, \mathbb{P})$. For a sub-$\sigma$-algebra $\mathcal{G}^0$ of $\mathcal{F}^0$, we denote the sub-$\sigma$-algebra $\mathcal{G}^0\otimes \{\varnothing,\Omega^1\}$ of $\mathcal{F}$  by $\mathcal{G}^0$ with a slight abuse of notation.

We now give some definitions  with respect to Lion's derivative.
\begin{De}\label{D1}
\rm
 The operator $\varphi: \mathcal{P}_{2}(\mathbb{R}^d) \rightarrow \mathbb{R}$ is called L-differential at $\mu \in \mathcal{P}_{2}(\mathbb{R}^d)$ if there exists a random variable $X \in \mathbb{L}_2(\mathbb{R}^d)$ such that 
 $\mu=\mathcal{L}(X)$ and the lifted function 
 $\tilde{\varphi}(X):=\varphi(\mathcal{L}(X))$ is $\mathrm{Fr\mathrm{\acute{e}}chet}$ differentiable at $X$.
 \end{De}
    The $\mathrm{Fr\mathrm{\acute{e}}chet}$ derivative  $\mathcal{D}\tilde{\varphi}(X)$, when seen as an element of $\mathbb{L}_2(\mathbb{R}^d)$
 via the Riesz-type representation theorem,  can be represented as $\mathcal{D}\tilde{\varphi}(X)=\partial_\mu\varphi(\mathcal{L}(X))(X)$,
 where $\partial_\mu\varphi(\mathcal{L}(X)): \mathbb{R}^d\ni y\rightarrow \partial_\mu\varphi(\mathcal{L}(X))(y)\in \mathbb{R}^d$ is mean square integral with respect to $\mathcal{L}(X)$ (see, e.g., \cite{CCD2022}).

%
\begin{De}
\rm
 Denote by $\tilde{C}^{2,2,1}:=\tilde{C}^{2,2,1}(\mathbb{R}^d\times\mathcal{P}_2(\mathbb{R}^d )\times\mathbb{R}_+;\mathbb{R})$ the family of all operators $\phi=\phi(x,\mu,t): 
\mathbb{R}^d \times \mathcal{P}_2(\mathbb{R}^d)\times\mathbb{R}_+  \rightarrow\mathbb{R}$ satisfying that
\begin{itemize}
\item[(i)]$\phi$ is twice continuously differentiable at $x$ and is continuously differentiable at $t$;
\item[(ii)] For any $(x, t)\in\mathbb{R}^d\times \mathbb{R}_+$, $\phi$ is twice L-differentiable at $\mu$ such that the L-derivative of $\phi$ at $\mu$: $\mathbb{R}^d \times \mathcal{P}_2(\mathbb{R}^d) \times\mathbb{R}_+ \times \mathbb{R}^d \ni(x,\mu, t,y) \mapsto \partial_\mu\phi(x,\mu,t)(y) \in \mathbb{R}^d$, the partial derivative of $\partial_\mu\phi(x,\mu,t)(y)$ at $x$:     ${\mathbb{R}^d \times }\mathcal{P}_2(\mathbb{R}^d) \times \mathbb{R}_+\times \mathbb{R}^d \ni(x,\mu,t, y) \mapsto \partial_x\partial_\mu \phi(x,\mu,t)(y) \in \mathbb{R}^d$, the derivative of $\partial_\mu\phi(x,\mu,t)(y)$ at $y$: $\mathbb{R}^d \times \mathcal{P}_2(\mathbb{R}^d) \times\mathbb{R}_+ \times \mathbb{R}^d \ni(x,\mu, t,y) \mapsto \partial_y\partial_\mu \phi(x,\mu,t)(y) \in \mathbb{R}^d$   and the L-derivative of $\partial_\mu\phi(x,\mu,t)(y)$ at $\mu$: $\mathbb{R}^d \times \mathcal{P}_2(\mathbb{R}^d) \times\mathbb{R}_+ \times \mathbb{R}^d \times \mathbb{R}^d \ni(x,\mu, t,y, z) \mapsto$ $\partial_\mu^2 \phi(x,\mu,t)(y, z){ \in \mathbb{R}^{d\times d}}$ all have the versions which are locally bounded and continuous at any points $(x,\mu,t,y)$ {with} 
    $y\in \mathrm{Supp}(\mu)$ and $(x,\mu,t,y,z)$ {with} 
    $y,z\in \mathrm{Supp}(\mu)$; 
  \item[(iii)]For any compact set $\mathcal{K} \subset \mathbb{R}^d \times \mathcal{P}_{2}(\mathbb{R}^d)\times\mathbb{R}_+$,
\begin{align*}
\sup _{(x,\mu,t) \in \mathcal{K}}\Big[&\int_{\mathbb{R}^d}|\partial_\mu \phi(x,\mu,t)(y)|^2 { \mu(\mathrm{d} y )}
+\int_{\mathbb{R}^d}|\partial_y \partial_\mu \phi(x,\!\mu,t)(y)|^2 \mu(\mathrm{d} y)
\\&
+\int_{\mathbb{R}^d} \int_{\mathbb{R}^d}|\partial_\mu^2 \phi(x,\mu,t)(y, z)|^2 \mu(\mathrm{d} y) \mu(\mathrm{d} z)
\\&+\int_{\mathbb{R}^d} |\partial_x \partial_\mu \phi(x,\mu,t)(y)|^2\mu(\mathrm{d} y) \Big]<+\infty .
\end{align*}
\end{itemize}
\end{De}

\section{The well-posedness and the It\^o formula}\label{Section6}
This paper focuses on the design of the delay feedback control function $u(t,X(t-\tau))$, where $\tau$ represents the response lag between the state observation and the system's response to the control, for the MV-SDEs with common noise \eqref{1} such that the controlled system
\eqref{6.1*} becomes bounded in infinite time horizon and further exponentially stable in the mean square.

Here, we use the simple form of the feedback control function with $u(t,x)=-\alpha x$ for $(t ,x)\in\R_+\times\R^d$, where $\alpha$ is a positive constant. Thus the controlled system \eqref{6.1*} becomes
\begin{align}\label{6.1}
\mathrm{d} X(t)
&=\big(f(X(t) ,\mathcal{L}^1(X(t)))-\alpha X(t-\tau)\big)
     \mathrm{d} t
+g(X(t),\mathcal{L}^1(X(t))) \mathrm{d} W_t
\nonumber\\&\quad
+g^0(X(t), \mathcal{L}^1(X(t))) \mathrm{d} W_t^0.
\end{align}
Suppose the initial data of the underlying system described by the MV-SDEs with common noise \eqref{6.1} is
\begin{align}\label{id}
\xi(\theta)=X(0), \quad -\tau\leq \theta\leq 0,
\end{align}
where the random variable $X(0)$ on $(\Omega^0,\mathcal{F}^0,\mathbb{P}^0)$ satisfies for $p\geq4$,
\begin{align}\label{2.8}
\mathbb{E}|X(0)|^p<\infty.
\end{align}
 In addition, we propose a hypothesis for the coefficients of the MV-SDEs with common noise \eqref{1}.
 \begin{assp}\label{A1}
 \rm There exists a positive constant $L$ such that for any $x,y\in\mathbb{R}^d$, $\mu,\nu\in\mathcal{P}_{2}(\mathbb{R}^d)$,
\begin{align*}
&|f(x, \mu)-f(y, \nu)|\vee|g(x, \mu)-g(y, \nu)|\vee|g^0(x, \mu)-g^0(y, \nu)|
\\&\leq L(|x-y|
+\mathbb{W}_2(\mu,\nu)).
\end{align*}
\end{assp}
By Assumption \ref{A1}, one can find a positive constant $A$ such that for any $x\in\mathbb{R}^d$, $\mu\in\mathcal{P}_{2}(\mathbb{R}^d)$,
\begin{align}\label{line}
\begin{aligned}
&|f(x, \mu)|\vee|g(x, \mu)|\vee|g^0(x, \mu)|
\leq A(1+|x|
+\mathbb{W}_2(\mu,\delta_0)).
\end{aligned}
\end{align}
Now we prepare the regularity of the solutions of the MV-SDEs with common noise \eqref{1} and the controlled system \eqref{6.1}.
\begin{lem}(see \cite[Vol-II, Proposition 2.8]{CD2018})
Under Assumption \ref{A1}, the MV-SDEs with common noise \eqref{1} with initial value $X(0)$ has a unique global solution $X(t)$ on $[0,\infty)$.
\end{lem}
\begin{thm}\label{pth}
Under Assumption \ref{A1}, the controlled system \eqref{6.1} with initial data \eqref{id} has a unique global solution $X(t)$ on $[0,\infty)$. 
\end{thm}
{\bf Proof.} We prove the existence and uniqueness of the solution of the controlled system \eqref{6.1} by the Banach fixed point theorem. Define the space $\mathbb{S}^{2, d}_{k,T}:=\mathbb{S}^{2, d}([(k-1)(T\wedge \tau),k(T\wedge \tau)])$ of all $\mathbb{R}^d$-valued $\mathcal{F}$-progressively measurable processes $X(t)$ on $t\in[(k-1)(T\wedge \tau),k(T\wedge \tau)]$ satisfying
$
\mathbb{E}[\sup_{(k-1)(T\wedge \tau) \leq t \leq k(T\wedge \tau)}|X(t)|^2]<\infty,
$
and we equip it with the norm
$
\|X\|_S^2=\mathbb{E}[\sup _{(k-1)(T\wedge \tau) \leq t \leq k(T\wedge \tau) }|X(t)|^2]<\infty .
$
The space $(\mathbb{S}^{2, d}_{k,T},\|\cdot\|_{\mathbb{S}})$ is a Banach space.
Consider $t\in[0,\tau]$. The system \eqref{6.1} becomes \begin{align}\label{s1}
\mathrm{d} X(t)
&=\big(f(X(t) ,\mathcal{L}^1(X(t)))-\alpha \xi(t-\tau)\big)
     \mathrm{d} t
+g(X(t),\mathcal{L}^1(X(t))) \mathrm{d} W_t
\nonumber\\&\quad
+g^0(X(t), \mathcal{L}^1(X(t))) \mathrm{d} W_t^0,
\end{align}
with initial value $X(0)$.  Furthermore, for all ${X} \in \mathbb{S}^{2, d}_{1,T}$, we can find a $\mathcal{P}_2(\mathbb{R}^d)$-valued and $\mathcal{F}^0$-adapted version of $(\mathcal{L}^1(X(t)))_{0 \leq t \leq T\wedge\tau}$ with continuous paths. Then we define for $0 \leq t \leq T\wedge\tau$,
\begin{align}\label{Y}
\begin{aligned}
Y(t)&=X(0) +\int_0^t \big(f( X(s), \mathcal{L}^1(X(s)))-\alpha \xi(t-\tau)\big) d s \\
&\quad +\int_0^t g(X(s), \mathcal{L}^1(X(s))) d W_s+\int_0^t g^0( X(s), \mathcal{L}^1(X(s))) d W_s^0.
\end{aligned}
\end{align}
Using the elementary inequality, the H\"older inequality, the Burkholder-Davis-Gundy inequality and \eqref{id}, we arrive at for any $T>0$,
 \begin{align*}
& \mathbb{E}\Big[ \sup _{0 \leq s \leq T\wedge\tau}|Y(s)|^2\Big]
 \\&\leq 4 \mathbb{E}|X(0)|^2+4 \mathbb{E} \left[\sup _{0 \leq s \leq T\wedge\tau}\Big|\int_0^s\big( f(X(r), \mathcal{L}^1(X(r)))-\alpha \xi(r-\tau)\big)\mathrm{d} r\Big|^2 \right]
\\
& \quad+4 \mathbb{E}\left[\sup _{0 \leq s \leq T\wedge\tau}\Big|\int_0^s g(X(r), \mathcal{L}^1(X(r)))\mathrm{d}W_r\Big|^2 \right]
 \\
& \quad+ 4\mathbb{E}\left[\sup _{0 \leq s \leq T\wedge\tau}\Big|\int_0^s g^0(X(r), \mathcal{L}^1(X(r)))\mathrm{d} W^0_r\Big|^2 \right]
\\
&\leq (4+8\alpha^2 T^2) \mathbb{E}|X(0)|^2+8 T\mathbb{E}\int_0^{T\wedge\tau}\big| f(X(s), \mathcal{L}^1(X(s)))\big|^2\mathrm{d}s
\\
& \quad+16\mathbb{E}\int_0^{T\wedge\tau} |g(X(s), \mathcal{L}^1(X(s)))|^2\mathrm{d} s
+ 16 \mathbb{E}\int_0^{T\wedge\tau} |g^0(X(s), \mathcal{L}^1(X(s)))|^2\mathrm{d}s.
\end{align*}
Noticing that
\begin{align}\label{ex}
\mathbb{E}(\mathbb{W}_2^2(\mathcal{L}^1(X(t)),\delta_0))\leq\mathbb{E}[\mathbb{E}^1|X(t)|^2]=\mathbb{E}|X(t)|^2
\end{align}
By \eqref{line}, the Fubini Theorem  and \eqref{ex}, we get that
 \begin{align*}
\begin{aligned}
\mathbb{E}\Big[ \sup _{0 \leq s \leq {T\wedge\tau}}|Y(s)|^2\Big]
&\leq 24(4+ T)A^2 T+(4+8\alpha^2 T^2) \mathbb{E}|X(0)|^2
\\&\quad+48( 4+ T)A^2\int_0^{T+\tau}\mathbb{E}|X(s)|^2\mathrm{d}s.
\end{aligned}
\end{align*}
Taking advantage of the Gronwall formula and \eqref{2.8}, we derive
\begin{align}\label{2.9}
\begin{aligned}
&\mathbb{E} \Big[\sup_{0 \leq t \leq {T\wedge\tau}}|X(t)|^2\Big] \leq \big(24(4 +T)A^2T+(4+8\alpha^2 T^2) \mathbb{E}|X(0)|^2\big) e^{48(4+T)A^2 T}<\infty,
\end{aligned}
\end{align}
which implies that $Y\in\mathbb{S}^{2,d}_{1,T}$. Define a mapping $\Phi:\mathbb{S}^{2,d}_{1,T}\ni X\rightarrow Y\in \mathbb{S}^{2,d}_{1,T}$. It is easy to see that the fixed point of $\Phi$ is the solution of the system \eqref{6.1} on $[0,T\wedge \tau]$. We proceed to find the fixed point of $\Phi$.  For processes $X,X^{\prime}\in\mathbb{S}^{2,d}_{1,T}$ with $X(0)=X^{\prime}(0)$, denote by $Y,Y^{\prime}$ the corresponding process defined by \eqref{Y}. Then $Y,Y^{\prime}\in\mathbb{S}^{2,d}_{1,T}$. The H\"older inequality, the Burkholder-Davis-Gundy inequality and \eqref{id} yields that for any $T>0$,
\begin{align*}
&\mathbb{E}\Big[\sup_{0 \leq t \leq {T\wedge\tau}}|Y(t)-Y^{\prime}(t)|^2\Big]
\\&\leq 3 \mathbb{E} \left[\sup _{0 \leq s \leq T\wedge\tau}\Big|\int_0^s\big( f(X^{\prime}(r), \mathcal{L}^1(X^{\prime}(r)))-f(X^{\prime}(r), \mathcal{L}^1(X^{\prime}(r)))\big)\mathrm{d} r\Big|^2 \right]
\\
& \quad+3 \mathbb{E}\left[\sup _{0 \leq s \leq T\wedge\tau}\Big|\int_0^s g(X(r), \mathcal{L}^1(X(r)))-g(X^{\prime}(r), \mathcal{L}^1(X^{\prime}(r)))\mathrm{d} W_r\Big|^2 \right]
 \\
& \quad+ 3\mathbb{E}\left[\sup _{0 \leq s \leq T\wedge\tau}\Big|\int_0^s g^0(X(r), \mathcal{L}^1(X(r)))-g^0(X^{\prime}(r), \mathcal{L}^1(X^{\prime}(r)))\mathrm{d} W^0_r\Big|^2 \right]
\\&\leq 3 T\mathbb{E}\int_0^{T\wedge\tau}\big| f(X(s), \mathcal{L}^1(X(s)))- f(X^{\prime}(s), \mathcal{L}^1(X^{\prime}(s)))\big|^2\mathrm{d}s
\\
& \quad+12\mathbb{E}\int_0^{T\wedge\tau} |g(X(s), \mathcal{L}^1(X(s)))-g(X^{\prime}(s), \mathcal{L}^1(X^{\prime}(s)))|^2\mathrm{d} s
\\
& \quad+ 12 \mathbb{E}\int_0^{T\wedge\tau} |g^0(X(s), \mathcal{L}^1(X(s)))-g^0(X^{\prime}(s), \mathcal{L}^1(X^{\prime}(s)))|^2\mathrm{d}s.
\end{align*}
This together with Assumption \ref{A1}  shows that
\begin{align*}
&\mathbb{E}\Big[\sup_{0 \leq t \leq {T\wedge\tau}}|Y(t)-Y^{\prime}(t)|^2\Big]
\\&\leq 6(8+ T)L^2\mathbb{E}\int_0^{T\wedge\tau}| X(s)-X^{\prime}(s)|^2+ \mathbb{W}_2^2(\mathcal{L}^1(X(s))-\mathcal{L}^1(X^{\prime}(s)))\mathrm{d}s
\end{align*}
Notice $$\mathbb{E}\mathbb{W}_2^2(\mathcal{L}^1(X(s)),\mathcal{L}^1(X^{\prime}(s)))\leq \mathbb{E}(\mathbb{E}^1|X(s)-X^{\prime}(s)|^2)=\mathbb{E}|X(s)-X^{\prime}(s)|^2.$$
Then we derive that
\begin{align}\label{contr}
\mathbb{E}\Big[\sup_{0 \leq t \leq {T\wedge\tau}}|Y(t)-Y^{\prime}(t)|^2\Big]
&\leq 12(8+ T)L^2T\mathbb{E}\Big[\sup_{0 \leq t \leq {T\wedge\tau}}| X(s)-X^{\prime}(s)|^2\Big].
\end{align}
Then one can easily find a $T_1$ such that $12(8+ T_1)L^2T_1<1$. Then the mapping $\Phi$ owns the unique fixed point on $\mathbb{S}^{2,d}_{1,T_1}$. If $T_1\geq\tau$, this implies that the system \eqref{6.1} exists a unique solution on $[0,\tau]$ satisfying $\mathbb{E}\big[\sup_{0\leq t\leq \tau}|X(t)|^2\big]<\infty$. Repeating the procedure on $[\tau ,2\tau], [2\tau, 3\tau],~\cdots$, we can obtain the unique solution $X(t)$ of the system \eqref{6.1} on $[0,\infty)$. If $T_1<\tau$, then we obtain that the system \eqref{6.1} exists a unique solution on $[0,T_1]$ satisfying $\mathbb{E}\big[\sup_{0\leq t\leq T_1}|X(t)|^2\big]<\infty$. Let $k=\lceil\tau/T_1\rceil$. Finding the fixed point of $\Phi$ on $\mathbb{S}^{2, d}_{2,T_1}$, we can derive the system \eqref{6.1} exists a unique solution on $[T_1,2T_1]$ satisfying $\mathbb{E}\big[\sup_{T_1\leq t\leq 2T_1}|X(t)|^2\big]<\infty$. Repeating the above procedure on $\mathbb{S}^{2, d}_{3,T_1},\cdots, \mathbb{S}^{2, d}_{k,T_1}$, one can derive the unique solution $X(t)$ of the system \eqref{6.1} on $[0,\tau]$. Repeating the argument above on $[\tau,2\tau], [2\tau, 3\tau],\cdots$, we can obtain the unique solution $X(t)$ of the system \eqref{6.1} on $[0,\infty)$. The proof is completed.\qed

We are now in a position to prove the It\^o formula with respect to both state and measure, the essential tool for the asymptotic analysis, for the controlled system \eqref{6.1}.
\begin{thm}\label{l*}
Let Assumption \ref{A1} hold and $V(x,\mu,t)\in \tilde{C}^{2,2,1}$. Then for the solution $X(t)$ of the controlled system \eqref{6.1}, it holds that 
\begin{align*}
&\mathbb{E}V(X(t),\mathcal{L}^1(X(t)),t)
\\&=\mathbb{E}V(X(0),\mathcal{L}^1(X(0)),0)
+\mathbb{E}\int_0^tLV(X(s),\mathcal{L}^1(X(s)),s)\mathrm{d} s,\quad  t\geq0,
\end{align*}
where
\begin{align}\label{LV}
&LV(X(s),\mathcal{L}^1(X(s)),s)
\nonumber\\&= \partial_t V(X(s),\mathcal{L}^1(X(s)),s)
+\partial_x V(X(s),\mathcal{L}^1(X(s)),s)\big(f(X(s) ,\mathcal{L}^1(X(s)))-\alpha X(s-\tau)\big)
\nonumber\\&+\frac{1}{2} \operatorname{tr}\big(\partial_{x x}V(X(s),\mathcal{L}^1(X(s)),s)g(X(s),\mathcal{L}^1(X(s)))
(g(X(s),\mathcal{L}^1(X(s))))^{\prime}\big)
\nonumber\\&+\frac{1}{2} \operatorname{tr}\big(\partial_{x x}V(X(s),\mathcal{L}^1(X(s)),s)g^0(X(s),\mathcal{L}^1(X(s)))
(g^0(X(s),\mathcal{L}^1(X(s))))^{\prime}\big)
\nonumber\\&+\mathbb{E}_1^1\Big\{\big(\partial_\mu V(X(s),\mathcal{L}^1(X(s)),s)(X_1(s))\big)^{\prime} \big(f(X_1(s),\mathcal{L}^1(X_1(s)))-\alpha X_1(s-\tau)\big)
\nonumber\\&+ \operatorname{tr}\big(\partial_x \partial_\mu V(X(s),\mathcal{L}^1(X(s)),s) (X_1(s)) g^0(X(s),\mathcal{L}^1(X(s)))(g^0(X_1(s),\mathcal{L}^1(X_1(s))))^{\prime}\big)
\nonumber\\&+\frac{1}{2} \operatorname{tr}\big(\partial_y \partial_\mu V(X(s),\mathcal{L}^1(X(s)),s)(X_1(s)) g(X_1(s),\mathcal{L}^1(X_1(s)))(g(X_1(s),\mathcal{L}^1(X_1(s))))^{\prime}\big)
\nonumber\\&+ \frac{1}{2} \operatorname{tr}\big(\partial_y \partial_\mu V(X(s),\mathcal{L}^1(X(s)),s)(X_1(s)) g^0(X_1(s),\mathcal{L}^1(X_1(s)))(g^0(X_1(s),\mathcal{L}^1(X_1(s))))^{\prime}\big)
\nonumber\\
&+\frac{1}{2} \mathbb{E}_2^1\big\{\operatorname{tr}\big(\partial_\mu^2 V(X(s),\mathcal{L}^1(X(s)),s)(X_1(s), {X_2}(s))g^0(X_1(s),\mathcal{L}^1(X_1(s)))
\nonumber\\&\quad\quad\quad\quad\quad\quad\quad\quad\quad\quad\quad\quad\times(g^0(X_2(s),\mathcal{L}^1(X_2(s))))^{\prime}\big)\big\}\Big\}
,\end{align}
$X_1(t)$ and $X_2(t)$ are copies of $X(t)$ on the spaces $\Omega_1$ and $\Omega_2$, respectively.
\end{thm}
{\bf Proof.}
 Let $X(t)$ be the unique solution of the controlled system \eqref{6.1}. Using the elementary inequality, the H\"older inequality and the Burkholder-Davis-Gundy inequality, we arrive at
 \begin{align*}
& \mathbb{E}\Big[ \sup _{0 \leq s \leq t}|X(s)|^2\Big] \\&\leq 4 \mathbb{E}|X(0)|^2+4 \mathbb{E} \left[\sup _{0 \leq s \leq t}\Big|\int_0^s\big( f(X(r), \mathcal{L}^1(X(r)))-\alpha X(r-\tau)\big)\mathrm{d} r\Big|^2 \right]
\\
& \quad+4 \mathbb{E}\left[\sup _{0 \leq s \leq t}\Big|\int_0^s g(X(r), \mathcal{L}^1(X(r)))\mathrm{d} B_r\Big|^2 \right]
 \\
& \quad+ 4\mathbb{E}\left[\sup _{0 \leq s \leq t}\Big|\int_0^s g^0(X(r), \mathcal{L}^1(X(r)))\mathrm{d} B^0_r\Big|^2 \right]
\\
&\leq 4 \mathbb{E}|X(0)|^2+4 T\mathbb{E}\int_0^t\big| f(X(s), \mathcal{L}^1(X(s)))-\alpha X(s-\tau)\big|^2\mathrm{d}s
\\
& \quad+16\mathbb{E}\int_0^t |g(X(s), \mathcal{L}^1(X(s)))|^2\mathrm{d} s
+ 16 \mathbb{E}\int_0^t |g^0(X(s), \mathcal{L}^1(X(s)))|^2\mathrm{d}s
\\&\leq 4 \mathbb{E}|X(0)|^2+8 T\mathbb{E}\int_0^t\big| f(X(s), \mathcal{L}^1(X(s)))\big|^2\mathrm{d}s
+8\alpha^2 T\mathbb{E}\int_0^t |X(s-\tau)|^2\mathrm{d} s
\\
& \quad+16\mathbb{E}\int_0^t |g(X(s), \mathcal{L}^1(X(s)))|^2\mathrm{d} s
+ 16 \mathbb{E}\int_0^t |g^0(X(s), \mathcal{L}^1(X(s)))|^2\mathrm{d}s.
\end{align*}
Integration by variable substitution gives
 \begin{align*}
\begin{aligned}
& \mathbb{E}\Big[ \sup _{0 \leq s \leq t}|X(s)|^2\Big]
\\&\leq 4 \mathbb{E}|X(0)|^2+8 T\mathbb{E}\int_0^t\big| f(X(s), \mathcal{L}^1(X(s)))\big|^2\mathrm{d}s
+8\alpha^2 T\mathbb{E}\int_{-\tau}^{t-\tau} |X(s)|^2\mathrm{d} s
\\
& \quad+16\mathbb{E}\int_0^t |g(X(s), \mathcal{L}^1(X(s)))|^2\mathrm{d} s
+ 16 \mathbb{E}\int_0^t |g^0(X(s), \mathcal{L}^1(X(s)))|^2\mathrm{d}s
\\&\leq 4(1 +2\alpha^2\tau T)\mathbb{E}|X(0)|^2+8 T\mathbb{E}\int_0^t\big| f(X(s), \mathcal{L}^1(X(s)))\big|^2\mathrm{d}s
\\
& \quad+8\alpha^2 T\mathbb{E}\int_{0}^{t} |X(s)|^2\mathrm{d} s
+16\mathbb{E}\int_0^t |g(X(s), \mathcal{L}^1(X(s)))|^2\mathrm{d} s
\\
& \quad+ 16 \mathbb{E}\int_0^t |g^0(X(s), \mathcal{L}^1(X(s)))|^2\mathrm{d}s.
\end{aligned}
\end{align*}
By \eqref{line}, the elementary inequality and the H\"older inequality, we get that
 \begin{align*}
\begin{aligned}
\mathbb{E}\Big[ \sup _{0 \leq s \leq t}|X(s)|^2\Big]
&\leq  4(1 +2\alpha^2\tau T)\mathbb{E}|X(0)|^2+8\alpha^2T\mathbb{E}\int_{0}^{t} |X(s)|^2\mathrm{d} s
\\
& \quad+8(4+ T)A^2\mathbb{E}\int_0^t( 1+|X(s)|+\mathbb{E}^1|X(s)| )^2\mathrm{d}s
\\&\leq  4(1 +2\alpha^2\tau T)\mathbb{E}|X(0)|^2+8\alpha^2T\mathbb{E}\int_{0}^{t} |X(s)|^2\mathrm{d} s
\\
& \quad+24(4+ T)A^2\mathbb{E}\int_0^t( 1+|X(s)|^2+\mathbb{E}^1|X(s)|^2 )\mathrm{d}s
.
\end{aligned}
\end{align*}
 The Fubini Theorem and \eqref{ex} yield that
  \begin{align*}
\mathbb{E}\Big[ \sup _{0 \leq s \leq t}|X(s)|^2\Big]
&\leq 24(4+ T)A^2T + 4(1 +2\alpha^2\tau T) \mathbb{E}|X(0)|^2
\\&\quad+8\big((24+6 T)A^2+\alpha^2 T\big)\int_0^t \mathbb{E}|X(s)|^2\mathrm{d}s
\\&\leq 24(4+ T)A^2T + 4(1 +2\alpha^2\tau T) \mathbb{E}|X(0)|^2
\\&\quad+8\big((24+6 T)A^2+\alpha^2 T\big)\int_0^t \mathbb{E}\Big[\sup_{0\leq r\leq s}|X(r)|^2\Big]\mathrm{d}s
.
\end{align*}
Taking advantage of the Gronwall formula and \eqref{2.8}, we derive for any $T>0$,
\begin{align}\label{2.9}
\begin{aligned}
&\mathbb{E} \Big[\sup _{0 \leq t \leq T}|X(t)|^2\Big] \\&\leq \big(24(4+ T)A^2T + 4(1 +2\alpha^2\tau T) \mathbb{E}|X(0)|^2\big) e^{8((24+6 T)A^2+\alpha^2 T) T}<\infty .
\end{aligned}
\end{align}
Using a similar method, one can obtain that for any $T>0$,
\begin{align}\label{2.10}
\mathbb{E} \Big[\sup _{0 \leq s \leq T}|X(s)|^4\Big] <\infty.
\end{align}
Using \eqref{line}, \eqref{2.9} and \eqref{2.10} yields that
\begin{align}\label{eq5.3*}
&\mathbb{E}\int_0^T\Big[|{f}(X(t),\mathcal{L}^1(X(t)))-\alpha X(t-\tau)|^2
+|{g}(X(t), \mathcal{L}^1(X(t)))|^4
\nonumber\\&\quad\quad\quad
+|{g}^0(X(t),\mathcal{L}^1(X(t)))|^4\Big]\mathrm{d}t<\infty.
\end{align}
Then by utilizing the It\^o formula with respect to measure (see, \cite[Vol-II, Theorem 4.17]{CD2018}) and taking expectation on both sides, the desired assertion holds.
\qed

 Next, we give some preparations for studying the asymptotic properties of system \eqref{6.1}. Let $X(\theta) :=X(0)$ for $-2\tau \leq \theta \leq -\tau$ and $\mathcal{L}^1(X(\theta)):=\mathcal{L}^1(X(0))$ for $-\tau \leq \theta \leq 0$. In order to estimate the derivation from the time delay in the mean square  $\mathbb{E}|X(t)-X(t-\tau)|^2$, we define an auxiliary function for $t\geq 0$,
 \begin{align}
I(t)=\int_{-\tau}^0\int_{t+r}^t&\big(\tau|f(X(s), \mathcal{L}^1(X(s)))-\alpha X(s-\tau)|^2+|g(X(s), \mathcal{L}^1(X(s)))|^2
\nonumber\\&+|g^0(X(s), \mathcal{L}^1(X(s))) |^2\big)\mathrm{d}s\mathrm{d}r.
\end{align}
It follows from a direct calculation that
\begin{align}\label{eq6.2}
  { \dot{I}(t)} =I_1(t)-I_2(t),
\end{align}
where
\begin{align*}
\begin{aligned}
I_1(t)&=\tau\big(\tau|f(X(t), \mathcal{L}^1(X(t)))-\alpha X(t-\tau)|^2
+|g(X(t), \mathcal{L}^1(X(t)))|^2\\&\quad+|g^0( X(t), \mathcal{L}^1(X(t))) |^2\big),
\\
I_2(t)&=\int_{t-\tau}^t\big(\tau|f(X(r) ,\mathcal{L}^1(X(r)))-\alpha X(r-\tau)|^2
+|g( X(r) ,\mathcal{L}^1(X(r)))|^2
\\&\quad+|g^0(X(r), \mathcal{L}^1(X(r))) |^2\big)\mathrm{d}r.
\end{aligned}
\end{align*}
By altering the integration  order, one derives that
\begin{align}\label{eq6.5}
I(t)\leq \tau I_2(t).
\end{align}
We proceed to estimate $\mathbb{E}|X(t)-X(t-\tau)|^2$. By employing the H\"older inequality and the It\^o isometry formula, one arrives at
\begin{align*}
&\mathbb{E}|X(t)-X(t-\tau)|^2\\&
\leq3\mathbb{E}\Big[\int_{t-\tau}^t \big(\tau|f( X(s),\mathcal{L}^1(X(s)))-\alpha X(s-\tau)|^2
\nonumber \\&  \quad
+|g(X(s) ,\mathcal{L}^1(X(s)))|^2
+|g^0( X(s), \mathcal{L}^1(X(s))|^2\big)\mathrm{d} s \Big]
\nonumber \\&
\leq 3\mathbb{E} I_2(t)
\end{align*}
for $t\geq\tau$, and \begin{align*}
&\mathbb{E}|X(t)-X(t-\tau)|^2=\mathbb{E}|X(t)-X(0)|^2\nonumber\\&
\leq3\mathbb{E}\Big[\int_{0}^t \big(\tau|f( X(s), \mathcal{L}^1(X(s)))-\alpha X(s-\tau)|^2
\nonumber \\&  \quad   +|g(X(s) ,\mathcal{L}^1(X(s))|^2+|g^0(X(s), \mathcal{L}^1(X(s))|^2\big)\mathrm{d} s \Big]
\nonumber\\&
\leq3\mathbb{E}\Big[\int_{t-\tau}^t \big(\tau|f(X(s), \mathcal{L}^1(X(s)))-\alpha X(s-\tau)|^2
\nonumber \\&  \quad +|g(X(s) ,\mathcal{L}^1(X(s)))|^2
+|g^0( X(s), \mathcal{L}^1(X(s))|^2\big)\mathrm{d} s \Big]
\nonumber \\&
\leq 3\mathbb{E}I_2(t)
\end{align*}
for $0\leq t\leq \tau$, which implies that
\begin{align}\label{eq6.6}
\mathbb{E}|X(t)-X(t-\tau)|^2\leq 3\mathbb{E} I_2(t).
\end{align}
\section{Boundedness control}
Since the  solution of system \eqref{1} may be unbounded on the infinite time horizon $[0,\infty)$ \cite{LMMY20}, it is necessary to design the delay feedback control function $u(t,X(t-\tau))$ such that the controlled system \eqref{6.1} becomes bounded in the mean square in $[0,\infty)$. 
To be precise, we give an useful estimation. Making use of Assumption \ref{A1}, one can find a positive constant $B$ such that for any $(x,\mu)\in\mathbb{R}^d \times\mathcal{P}_{2}(\mathbb{R}^d)$
\begin{align}\label{mono}
 &2\langle x,  f( x, \mu)\rangle+|g(x, \mu)|^2 +|g^0( x, \mu)|^2
\leq B(1+|x|^2+\mathbb{W}^2_2(\mu,\delta_0)).
\end{align}
 \begin{thm}\label{T6.1}
Let {\rm Assumption \ref{A1}} hold and $ \alpha>B$. If  $0<\tau<\tau ^*:=\tau_1\wedge \tau_2$, then the solution of system \eqref{6.1} has the property that
\begin{align}\label{6.4}
\sup_{0\leq t<\infty}\mathbb{E}|X(t)|^2< \infty,
\end{align}
where
$ \tau_1$ and $\tau_2$ are the positive solutions of equations
$$\phi(\tau_1)=\frac{1}{6},~~~~\psi(\tau_2)=\frac{( \alpha-B)^2 }{6\alpha^2},$$
 respectively. Here $\phi(\cdot)$ and $\psi(\cdot)$ are defined by \eqref{eq6.8*} and \eqref{eq6.12*}, respectively.
\end{thm}
{\bf Proof.} Fix $0<\tau<{\tau^*}$.
Consider
\begin{align}\label{eq6.3}
&V(X(t),\mathcal{L}^{1}(X(t)),t)=|X(t)|^2+\mathbb{E}^1|X(t)|^2+\sigma I(t),
\end{align}
where $I(t)$ is defined by \eqref{eq6.2} and $\sigma$ is a positive constant to be determined later.
In order to estimate $V({{X}}(t),\mathcal{L}^{1}(X(t)), t)$, we first analyze $I(t)$.
By the definition of $I_1(t)$, the elementary inequality, \eqref{line}  as well as  the H\"older inequality,
we have
\begin{align*}
I_1(t)
&\leq \tau\big(2\tau|f(X(t), \mathcal{L}^1(X(t)))|^2+2\alpha^2\tau |X(t-\tau)|^2
\\&\quad+|g(X(t), \mathcal{L}^1(X(t)))|^2+|g^0( X(t), \mathcal{L}^1(X(t)))|^2\big)
\\&\leq \tau\big(2\tau|f(X(t), \mathcal{L}^1(X(t)))|^2+4\alpha^2\tau |X(t)|^2+4\alpha^2\tau |X(t)-X(t-\tau)|^2
\\&\quad+|g(X(t), \mathcal{L}^1(X(t)))|^2+|g^0( X(t), \mathcal{L}^1(X(t)))|^2\big)
\\&\leq \tau\big((6 A^2\tau+6 A^2)(1+|X(t)|^2
+\mathbb{E}^1|X(t)|^2)
\\&\quad+4\alpha^2\tau |X(t)|^2+4\alpha^2\tau |X(t)-X(t-\tau)|^2\big)
\\
&\leq\tau\Big[\big(6A^2+(6 A^2+4\alpha^2)\tau\big)|X(t)|^2
+6 A^2(1+\tau)\mathbb{E}^1|X(t)|^2
\\&\quad+4\alpha^2\tau|X(t)-X(t-\tau)|^2
+6 A^2(1 +\tau)
\Big].
\end{align*}
This combined with \eqref{eq6.2} gives that
\begin{align}\label{eq6.8}
\dot{I}(t)&\leq\beta_1(\tau)|X(t)|^2
+\beta_2(\tau)\mathbb{E}^1|X(t)|^2
+\phi(\tau)|X(t)-X(t-\tau)|^2
\nonumber\\&\quad+\beta_2(\tau)-I_2(t),
\end{align}
where \begin{align}\label{eq6.8*}
  &\phi(\tau)=4\alpha^2\tau^2,\\
 \beta_1(\tau) =6A^2\tau +2(3A^2&+2\alpha^2) \tau^2, ~ \beta_2(\tau)=6A^2  (\tau +\tau^2). \nonumber
\end{align}
We now turn to estimate \begin{align}\label{U}
U(X(t),\mathcal{L}^1(X(t))):=|X(t)|^2+\mathbb{E}^1|X(t)|^2.
\end{align}
Applying Lemma \ref{l*} to $U(X(t),\mathcal{L}^1(X(t)))$ and noticing that $X(t)$ and $X_1(t)$ are identically distributed show that
\begin{align*}
&LU(X(t),\mathcal{L}^1(X(t)))
\\&=2\langle X(t),f(X(t) ,\mathcal{L}^1(X(t)))-\alpha X(t-\tau)\rangle
\\&\quad+|g(X(t),\mathcal{L}^1(X(t)))|^2+|g^0(X(t),\mathcal{L}^1(X(t)))|^2
\\&\quad+\mathbb{E}_1^1 \big\{2\langle X_1(t),f(X_1(t) ,\mathcal{L}^1(X_1(t)))-\alpha X_1(t-\tau)\rangle
\\&\quad+|g(X_1(t),\mathcal{L}^1(X_1(t)))|^2+|g^0(X_1(t),\mathcal{L}^1(X_1(t)))|^2\big\}
\\&=2\langle X(t),f(X(t),\mathcal{L}^1(X(t)))-\alpha X(t-\tau)\rangle
\\&\quad+|g(X(t),\mathcal{L}^1(X(t)))|^2+|g^0(X(t),\mathcal{L}^1(X(t)))|^2
\\&\quad+\mathbb{E}^1 \big\{2\langle X(t),f(X(t) ,\mathcal{L}^1(X(t)))-\alpha X(t-\tau)\rangle
\\&\quad+|g(X(t),\mathcal{L}^1(X(t)))|^2+|g^0(X(t),\mathcal{L}^1(X(t)))|^2\big\}
\end{align*}
Making use of \eqref{mono} as well as the elementary inequality, one obtains that
\begin{align}\label{eq3.5}
&LU(X(t),\mathcal{L}^1(X(t)))
\nonumber\\&\leq B(1+|X(t)|^2
+\mathbb{E}^1|X(t)|^2)
-2\alpha |X(t)|^2+ 2\alpha\langle X(t), X(t)-X(t-\tau)\rangle
\nonumber\\&\quad+\mathbb{E}^1 \big\{B(1+|X(t)|^2
+\mathbb{E}^1|X(t)|^2)
-2\alpha |X(t)|^2+ 2\alpha\langle X(t), X(t)-X(t-\tau)\rangle\big\}
\nonumber\\&\leq -\Big(2\alpha-B-\frac{\zeta}{2}\Big)| X(t)|^2
+\frac{2\alpha^2 }{\zeta}|X(t)-X(t-\tau)|^2
-\Big(2\alpha-3B-\frac{\zeta}{2}\Big)\mathbb{E}^1| X(t)|^2
\nonumber\\&\quad+\frac{2\alpha^2 }{\zeta}\mathbb{E}^1|X(t)-X(t-\tau)|^2
+ 2B,
\end{align}
for any $\zeta>0$. Combining \eqref{eq6.8} with \eqref{eq3.5} yields that
\begin{align}\label{3.6}
&LV(X(t),\mathcal{L}^{1}(X(t)),t)=LU(X(t),\mathcal{L}^{1}(X(t)))+\sigma\dot{I}(t)
\nonumber\\
&
\leq
-\Big(2\alpha-B-\frac{\zeta}{2}-\sigma\beta_1(\tau)\Big)| X(t)|^2+\Big(\frac{2\alpha^2 }{\zeta}+\sigma\phi(\tau)\Big)|X(t)-X(t-\tau)|^2
\nonumber\\&\quad
-\Big(2\alpha-3B-\frac{\zeta}{2}-\sigma\beta_2(\tau)\Big)\mathbb{E}^1| X(t)|^2+\frac{2\alpha^2 }{\zeta}\mathbb{E}^1|X(t)-X(t-\tau)|^2
\nonumber\\&\quad
-\sigma I_2(t)+\sigma\beta_2(\tau)+ 2B.
\end{align}
By Lemma \ref{l*}, \eqref{eq6.5} and \eqref{3.6}, 
one has
 \begin{align*}
&\mathbb{E}[e^{\gamma t}V(X(t),\mathcal{L}^{1}(X(t)),t)]
\\&=2\mathbb{E} |X(0)|^2+\sigma \mathbb{E} I(0)+\mathbb{E}\int_{0}^te^{\gamma s}(\gamma V(X(s),\mathcal{L}^{1}(X(s)),s)+LV(X(s),\mathcal{L}^{1}(X(s)),s))\mathrm{d}s
\\
&
\leq2\mathbb{E} |X(0)|^2+\sigma \mathbb{E} I(0)+\mathbb{E}\int_{0}^te^{\gamma s}\Big(
-\Big(2\alpha-B-\frac{\zeta}{2}-\sigma\beta_1(\tau)-\gamma\Big)| X(s)|^2
\nonumber\\&\quad
+\Big(\frac{2\alpha^2 }{\zeta}+\sigma\phi(\tau)\Big)|X(s)-X(s-\tau)|^2-\Big(2\alpha-3B-\frac{\zeta}{2}-\sigma\beta_2(\tau)-\gamma\Big)\mathbb{E}^1| X(s)|^2
\nonumber\\&\quad
+\frac{2\alpha^2 }{\zeta}\mathbb{E}^1|X(s)-X(s-\tau)|^2-(\sigma- \sigma\tau\gamma) I_2(s)+\sigma\beta_2(\tau)+ 2B\Big)\mathrm{d}s,
\end{align*}
for some $\gamma>0$. Utilizing the Fubini Theorem and \eqref{ex}, we get
\begin{align*}
&\mathbb{E}[e^{\gamma t}V(X(t),\mathcal{L}^{1}(X(t)),t)]
\\
&
\leq2\mathbb{E} |X(0)|^2+\sigma \mathbb{E} I(0)+\int_{0}^te^{\gamma s}\Big(
-\Big(2\alpha-B-\frac{\zeta}{2}-\sigma\beta_1(\tau)-\gamma\Big)
\mathbb{E}| X(s)|^2
\nonumber\\&\quad
+\Big(\frac{2\alpha^2 }{\zeta}+\sigma\phi(\tau)\Big)\mathbb{E}|X(s)-X(s-\tau)|^2
-\Big(2\alpha-3B-\frac{\zeta}{2}-\sigma\beta_2(\tau)-\gamma\Big)\mathbb{E}\big[\mathbb{E}^1| X(s)|^2\big]
\nonumber\\&\quad
+\frac{2\alpha^2 }{\zeta}\mathbb{E}\big[\mathbb{E}^1|X(s)-X(s-\tau)|^2\big]-(\sigma- \sigma\tau\gamma) I_2(s)+\sigma\beta_2(\tau)+ 2B\Big)\mathrm{d}s
\\
&
\leq2\mathbb{E} |X(0)|^2+\sigma \mathbb{E} I(0)+\int_{0}^te^{\gamma s}\Big(
-(4\alpha-4B
-{\zeta}- \sigma\psi (\tau)-2\gamma)
\mathbb{E}| X(s)|^2
\nonumber\\&\quad
+\Big(\frac{4\alpha^2 }{\zeta}+\sigma\phi(\tau)\Big)\mathbb{E}|X(s)-X(s-\tau)|^2
-(\sigma- \sigma\tau\gamma) I_2(s)+\sigma\beta_2(\tau)+ 2B\Big)\mathrm{d}s
\end{align*}
where \begin{align}\label{eq6.12*}
\psi (\tau)= 12A^2 \tau+ 4(3A ^2+\alpha^2)\tau^2.
\end{align}
It follows from \eqref{eq6.6} that
\begin{align*}
&\mathbb{E}[e^{\gamma t}V(X(t),\mathcal{L}^{1}(X(t)),t)]
\\&\leq2\mathbb{E} |X(0)|^2+\sigma \mathbb{E} I(0)
+\int_0^t e^{\gamma s}\big[-(4\alpha-4B
-{\zeta}- \sigma\psi (\tau)-2\gamma)\mathbb{E}|X(s)|^2
\\
&\quad
-\Big(\sigma(1-3\phi(\tau))-\frac{12\alpha^2 }{\zeta}-\sigma\tau\gamma\Big) \mathbb{E} I_2(s)
+\sigma\beta_2(\tau)+ 2B\big]\mathrm{d}s,
\end{align*}
Choose
\begin{align}\label{eq6.13*}
\zeta=2(\alpha-B), \quad \sigma= \frac{12\alpha^2}{\alpha-B}.
\end{align}
Due to the increasing property of $\psi (\tau)$ and $\phi(\tau)$ in $\tau>0$ as well as the value of $\zeta$ and $\sigma$, we have $4\alpha-4B
-{\zeta}-\sigma\psi (\tau)>0, \sigma(1-3\phi(\tau))-{6\alpha^2 }/{(\alpha-B)}>0.$ Then we can find a sufficient small $\gamma$ such that
\begin{align*}
4\alpha-4B
-{\zeta}-\sigma\psi (\tau)-2\gamma&>0,\quad
 \sigma(1-3\phi(\tau))-\frac{6\alpha^2 }{\alpha-B}-\sigma\tau \gamma>0.
\end{align*}
Then, this together with   $\sigma\beta_2(\tau)\leq \sigma\psi(\tau)/2 \leq  {\alpha-B}$  implies
\begin{align*}
&\mathbb{E}[e^{\gamma t}V(X(t),\mathcal{L}^{1}(X(t)),t)]
\leq 2\mathbb{E}|X(0)|^2+\sigma\mathbb{E}[I(0)]
+ \frac{ \alpha+B}{ \gamma} e^{\gamma t}.
\end{align*}
Dividing $e^{\gamma t}$ on both side and using the definition of $V({{X}}_t,{\mathcal{L}}^{1, X}_t, t)$, we get
\begin{align*}
\mathbb{E}|X(t)|^2
&\leq 2\mathbb{E}|X(0)|^2 e^{-\gamma t}+\sigma\mathbb{E}[I(0)]e^{-\gamma t}
  +  \frac{ \alpha+B}{ \gamma}.
\end{align*}
Therefore, the required assertion follows. \qed
\section{Stabilization}
This section pays attention to the stabilization problem of the MV-SDEs with common noise \eqref{1}. To guarantee existence of the trivial solution, we further suppose $f(0,\delta_0)=0$, $g(0,\delta_0)=0$, $g^0(0,\delta_0)=0$ in this section which together with Assumption \ref{A1} gives that
there exist positive constants $C$ and $D$ such that for any $(x,\mu)\in\mathbb{R}^d \times\mathcal{P}_{2}(\mathbb{R}^d)$,
\begin{align}\label{line1}
&|f(x, \mu)|\vee|g(x, \mu)|\vee|g^0(x, \mu)|
\leq C(|x|
+\mathbb{W}_2(\mu,\delta_0)),
\end{align}
and
\begin{align}\label{mon1}
 2\langle x,  f(x, \mu)\rangle+|g(x, \mu)|^2+|g^0(x, \mu)|^2
  \leq D(|x|^2+\mathbb{W}_2^2(\mu,\delta_0)).
\end{align}
\begin{thm}\label{T6.2}
Let {\rm Assumption \ref{A1}}  hold and $\alpha>D$.   Then for any  $0<\tau<\tau^{**}:=\tau_3\wedge \tau_4$, the solution of system \eqref{6.1} has the property that
\begin{align}\label{T62*}
\limsup_{t\rightarrow\infty}\frac{1}{t}\log(\mathbb{E}|X(t)|^2)\leq -\gamma,
\end{align}
where
\begin{align}\label{eq6.15*}
\gamma=\frac{1-6\phi(\tau)}{2\tau}\wedge\frac{(\alpha-D)^2-6\alpha^2\varphi(\tau)}{\alpha-D},
\end{align}
and
$ \tau_3$, $\tau_4$ are the positive solutions of  equations
\begin{align}\label{eq6.15**}
\phi(\tau_3)=\frac{1}{6},~~~\quad\varphi(\tau_4)=\frac{(\alpha-D)^2}{6\alpha^2},
\end{align}
respectively. Here $\phi(\cdot)$ and $\varphi(\cdot)$ are defined by \eqref{eq6.8*} and \eqref{eq6.21*}, respectively.
\end{thm}
{\bf Proof.} Fix $0<\tau<\tau^{**}$. Applying \eqref{line1} and the elementary inequality, we derive from the definition of ${I}_1(t)$ that
 \begin{align*}
I_1(t)
&\leq 4 \tau \big(C^2 +  (C^2+ \alpha^2) \tau \big)|X(t)|^2
+4 C^2 \tau (1+\tau) \mathbb{E}^1|X(t)|^2
\nonumber\\&
\quad+4\alpha^2\tau^2|X(t)-X(t-\tau)|^2.
\end{align*}
This together with \eqref{eq6.2} yields that
\begin{align}\label{eq6.16}
\dot{I}( t)&\leq \beta_3(\tau)|X(t)|^2+\beta_4(\tau)\mathbb{E}^1|X(t)|^2
+\phi(\tau)|X(t)-X(t-\tau)|^2
-I_2(t),
\end{align}
where $\phi(\tau)$ is defined in \eqref{eq6.8*}, and \begin{align*}
\begin{aligned}
&\beta_3(\tau)=4C^2\tau+4(C^2+\alpha^2)\tau^2,~ \beta_4(\tau)=4C^2(\tau+\tau^2).
\end{aligned}
\end{align*}
We now turn to estimate $
U(X(t),\mathcal{L}^1(X(t)))
$ defined in \eqref{U}.
Applying Lemma \ref{l*}, \eqref{mon1} and the elementary inequality, one obtains that
\begin{align}\label{eq4.5}
&LU(X(t),\mathcal{L}^1(X(t)))
\nonumber\\&\leq D(|X(t)|^2
+\mathbb{E}^1|X(t)|^2)
-2\alpha |X(t)|^2+ 2\alpha\langle X(t), X(t)-X(t-\tau)\rangle
\nonumber\\&\quad+\mathbb{E}^1 \big\{B(|X(t)|^2
+\mathbb{E}^1|X(t)|^2)
-2\alpha |X(t)|^2+ 2\alpha\langle X(t), X(t)-X(t-\tau)\rangle\big\}
\nonumber\\&\leq -\Big(2\alpha-D-\frac{\zeta}{2}\Big)| X(t)|^2
+\frac{2\alpha^2 }{\zeta}|X(t)-X(t-\tau)|^2
-\Big(2\alpha-3D-\frac{\zeta}{2}\Big)\mathbb{E}^1| X(t)|^2
\nonumber\\&\quad+\frac{2\alpha^2 }{\zeta}\mathbb{E}^1|X(t)-X(t-\tau)|^2
,
\end{align}
for any $\zeta>0$.
For $
V(X(t),\mathcal{L}^1(X(t)),t)$
 given by \eqref{eq6.3}, recalling \eqref{eq6.16} and \eqref{eq4.5}, one arrives at
\begin{align*}
&LV(X(t),\mathcal{L}^{1}(X(t)),t)
\nonumber\\
&
\leq
-\Big(2\alpha-D-\frac{\zeta}{2}-\sigma\beta_3(\tau)\Big)| X(t)|^2+\Big(\frac{2\alpha^2 }{\zeta}+\sigma\phi(\tau)\Big)|X(t)-X(t-\tau)|^2
\nonumber\\&\quad
-\Big(2\alpha-3D-\frac{\zeta}{2}-\sigma\beta_4(\tau)\Big)\mathbb{E}^1| X(t)|^2+\frac{2\alpha^2 }{\zeta}\mathbb{E}^1|X(t)-X(t-\tau)|^2
\nonumber\\&\quad
-\sigma I_2(t),
\end{align*}
where $\zeta,~\sigma$ are defined by \eqref{eq6.13*}.
This together with Lemma \ref{l*} and \eqref{eq6.5} implies that
 \begin{align*}
&\mathbb{E}[e^{\gamma t}V(X(t),\mathcal{L}^{1}(X(t)),t)]
\\&
\leq2\mathbb{E} |X(0)|^2+\sigma \mathbb{E} I(0)+\mathbb{E}\int_{0}^te^{\gamma s}\Big(
-\Big(2\alpha-D-\frac{\zeta}{2}-\sigma\beta_3(\tau)-\gamma\Big)| X(s)|^2
\nonumber\\&\quad
+\Big(\frac{2\alpha^2 }{\zeta}+\sigma\phi(\tau)\Big)|X(s)-X(s-\tau)|^2
-\Big(2\alpha-3D-\frac{\zeta}{2}-\sigma\beta_4(\tau)-\gamma\Big)\mathbb{E}^1| X(s)|^2
\nonumber\\&\quad
+\frac{2\alpha^2 }{\zeta}\mathbb{E}^1|X(s)-X(s-\tau)|^2-(\sigma- \sigma\tau\gamma) I_2(s)\Big)\mathrm{d}s,
\end{align*}
for some $\gamma>0$.
 By the Fubini Theorem, \eqref{ex} and \eqref{eq6.6}, we derive that
\begin{align}\label{eq6.21}
&\mathbb{E}[e^{\gamma t}V(X(t),\mathcal{L}^1(X(t)),t)]\nonumber
\\&
\leq2\mathbb{E} |X(0)|^2+\sigma \mathbb{E} I(0)+\mathbb{E}\int_{0}^te^{\gamma s}\Big(
-(4\alpha-4D-\zeta
-\sigma\varphi(\tau)-2\gamma)\mathbb{E}|X(s)|^2
\nonumber\\&\quad
+\Big(\frac{4\alpha^2 }{\zeta}+\sigma\phi(\tau)\Big)\mathbb{E}|X(s)-X(s-\tau)|^2
-(\sigma- \sigma\tau\gamma) I_2(s)\Big)\mathrm{d}s
\nonumber\\&\leq2\mathbb{E}|X(0)|^2+\sigma\mathbb{E}[I(0)]
+\int_0^t e^{\gamma s}\big[-(4\alpha-4D-\zeta
-\sigma\varphi(\tau)-2\gamma)\mathbb{E}|X(s)|^2
\nonumber\\&\quad
-\Big(\sigma(1-3\phi(\tau))-\frac{12\alpha^2 }{\zeta}-\sigma\tau\gamma\Big) \mathbb{E} I_2(s)\big]\mathrm{d}s.
\end{align}
where
\begin{align}\label{eq6.21*}
\varphi(\tau)=8C^2\tau+4(2C^2+\alpha^2)\tau^2.
\end{align}
It follows from the increasing  property of $\phi(\tau)$ and $\varphi(\tau)$ that $4\alpha-4D-\zeta
-\sigma\varphi(\tau)>0,
~\sigma(1-3\phi(\tau))-12\alpha^2 /\zeta>0.$
One observes from \eqref{eq6.15*} that
\begin{align*}
&4\alpha-4D-\zeta
-\sigma\varphi(\tau)-2\gamma\geq0,
\quad\sigma(1-3\phi(\tau))-\frac{12\alpha^2 }{\zeta}-\sigma\tau\gamma\geq0.\end{align*}
 Dividing $e^{\gamma t}$ on both sides of \eqref{eq6.21}  and recalling the definition of $V(X(t),\mathcal{L}^1(X(t)),t)$ arrives at
\begin{align*}
\mathbb{E}|X(t)|^2
\leq (2\mathbb{E}|X(0)|^2+\sigma\mathbb{E}[I(0)])e^{-\gamma t}.
\end{align*}
The proof is therefore complete.\qed

We provide an example to illustrate the result of the stabilization.
\begin{figure}[htbp]
\centering
\subfigure{\includegraphics[height=5cm,width=6.5cm]{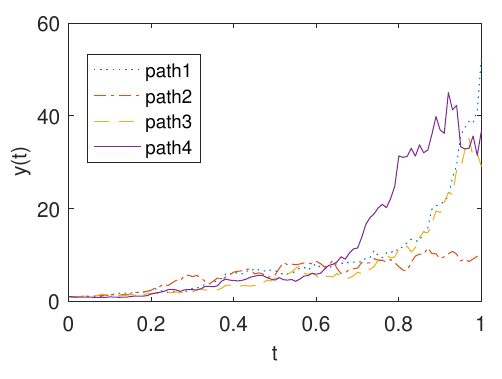}}
\quad\quad
\subfigure{\includegraphics[height=5cm,width=7cm]{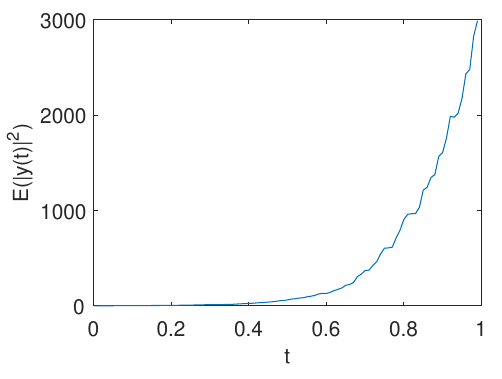}}
\vspace{-1em}
\caption{Four sample paths of $y(t)$ and the sample mean of $|y(t)|^2$  of system \eqref{eq7.2} for $t\in[0,1]$ with  $50$ sample points and step size $\Delta=0.01$.}
\end{figure}
\begin{expl}
\rm
Consider the unstable scalar system 
\begin{align}\label{eq7.2}
\mathrm{d} y(t)&=\big(3y(t)+\mathbb{E}^1y(t)\big) \mathrm{d} t+\big( y(t)+\mathbb{E}^1y(t)\big) \mathrm{d} W_t
\nonumber\\&\quad
+\big( y(t)+\mathbb{E}^1y(t)\big) \mathrm{d} W^0_t,
\end{align}
with initial data $y(0)=1$.
 Figure 1 depicts four sample paths of $y(t)$ and the sample mean of $|y(t)|^2$  of  \eqref{eq7.2} for $t\in[0,1]$ with   $50$ sample points and step size $\Delta=0.01$. Obviously, the coefficients are  globally Lipschitz continuous. It is easy to see Assumption \ref{A1} holds with $L=3$. We compute that
\begin{align*}
&|f(x,\mu)|\vee|g(x,\mu)|\vee|g^0(x,\mu)|
\leq 3(|x|+\mathbb{W}_2(\mu,\delta_0)),\\
&2\langle x,  f( x, \mu)\rangle+|g(x, \mu)|^2+|g^0( x, \mu)|^2
\leq 11(|x|^2+\mathbb{W}_2^2(\mu,\delta_0)).
\end{align*}
That is $C=3$ and $D=11$.
Choose $\alpha=22>D.$
The controlled system becomes
\begin{align}\label{eq7.3}
\mathrm{d} y(t)&=\big(3 y(t)-22 y(t-\tau)+\mathbb{E}^1y(t)\big) \mathrm{d} t
+\big( y(t)+\mathbb{E}^1y(t)\big) \mathrm{d} W_t
\nonumber\\&\quad
+\big(y(t)+\mathbb{E}^1y(t)\big) \mathrm{d} W^0_t,
\end{align}
where the initial data  $y(\theta)=1$, $-\tau\leq \theta\leq 0$.
By \eqref{eq6.15**} we obtain $\tau^{**}\approx5.696\times 10^{-4}$. Choose $\tau=5\times 10^{-4}$.  By virtue of Theorem \ref{T6.2}, the solution $y(t)$ of the controlled system \eqref{eq7.3} satisfies
$$
\limsup_{t\rightarrow\infty}\frac{1}{t}\log[\mathbb{E}|y(t)|^2]\leq-  1.363.
$$
 Figure 2 depicts four sample paths of the solution $y(t)$ and the sample mean of $|y(t)|^2$  of the controlled system \eqref{eq7.3} for $t\in[0,0.6]$  with 50 sample points and step size $\Delta=5\times 10^{-4}$.
 \begin{figure}[htbp]\label{f2}
\centering
\subfigure{\includegraphics[height=5cm,width=7cm]{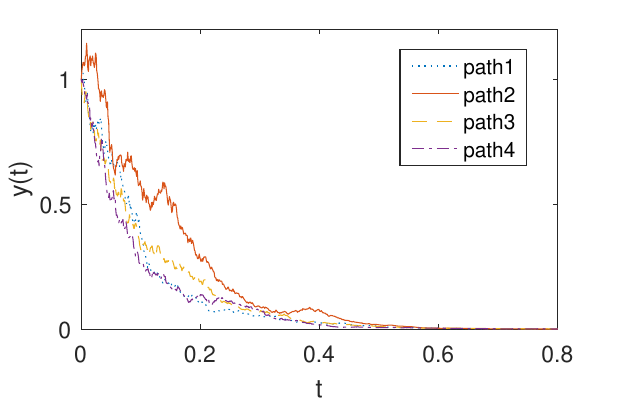}}
\quad\quad
\subfigure{\includegraphics[height=5cm,width=7cm]{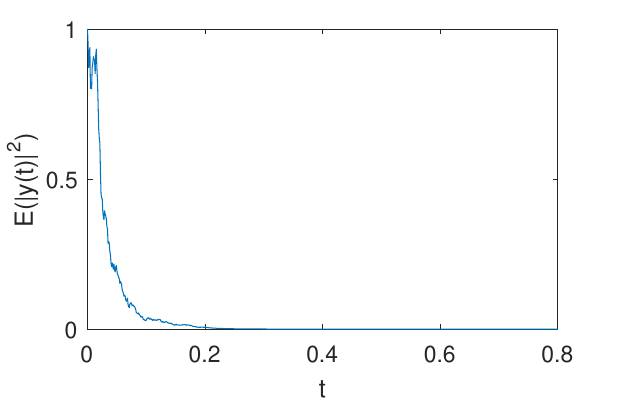}}
\vspace{-1em}
\caption{Four sample paths of $y(t)$ and the sample mean of $|y(t)|^2$ of the controlled system \eqref{eq7.3} for $t\in[0,0.6]$ with 50 sample points and step size $\Delta=5\times 10^{-4}$.}
\end{figure}
\end{expl}
\section{Conclusions}
In this paper, we design a delay feedback control term that is solely contingent on the state, which is more realistic and straightforward to implement. We prove the existence and uniqueness of the global solution of the controlled system. The It\^o formula with respect to both state and measure of the controlled system is derived. Then the control mechanism designed here ensures that the controlled system exhibits asymptotic boundedness and mean-square exponential stability. Finally, the permissible range of the delay $\tau$ is quantified.


\end{document}